\documentclass[preprint,11pt,authoryear]{elsarticle}
\usepackage{amssymb}

\usepackage[colorinlistoftodos]{todonotes}

\usepackage{latexsym}
\usepackage{amsfonts}
\usepackage{textcomp}
\usepackage{graphicx}
\usepackage{amssymb, amsmath}

\usepackage{color}

\newcommand{\A}{{\mathcal A}}

\newcommand{\E}{\mathbb{E}}
\renewcommand{\P}{\mathbb{P}}

\newcommand{\R}{\mathbb{R}}

\newcommand{\Z}{\mathbb{Z}}
\newcommand{\N}{\mathbb{N}}

\newcommand{\FF}{\mathcal F}

\newtheorem{theorem}{Theorem}
\newtheorem{corollary}{Corollary}

\newtheorem{lemma}{Lemma}
\newtheorem{proposition}{Proposition}
\newtheorem{remarque}{Remark}

\begin{document}

\begin{frontmatter}

\title{Convergence of generalized urn models to non-equilibrium attractors}

\author{Mathieu Faure}
\address{Aix-Marseille University (Aix-Marseille School of Economics)\\
CNRS and EHESS \\
Centre de la vieille Charit\'{e}\\ 2 rue de la Charit\'{e}, 13002 Marseille, France}
\author{Sebastian J. Schreiber}
\address{Department of Evolution and Ecology\\ University of California, Davis, California USA 95616}

\begin{abstract} Generalized Polya urn models have been used to model the establishment dynamics of a small founding population consisting of $k$ different genotypes or strategies. As population sizes get large, these population processes are well-approximated by a mean limit ordinary differential equation whose state space is the $k$ simplex. We prove that if this mean limit ODE has an attractor at which the temporal averages of the population growth rate is positive, then there is a positive probability of the population not going extinct (i.e. growing without bound) and its distribution converging to the attractor. Conversely, when the temporal averages of the population growth rate is negative along this attractor, the population distribution does not converge to the attractor. For the stochastic analog of the replicator equations which can exhibit non-equilibrium dynamics, we show that verifying the conditions for convergence and non-convergence reduces to a simple algebraic problem. We also apply these results to selection-mutation dynamics to illustrate convergence to periodic solutions of these population genetics models with positive probability. \vskip 0.1in

\end{abstract}

\end{frontmatter}

\section{Introduction} \label{sec:intro}
Biological invasions, where a species is introduced in a novel habitat, are occurring repeatedly throughout the world and often start with small founding population. Whether or not these founding populations establish or go extinct in their new environment depends on a diversity of factors including  local environmental conditions,  demographic stochasticity, genetic diversity of the founding population, and nonlinear feedbacks between individuals in the founding population. One commonly used approach to understanding the roles of the first two factors is modeling the dynamics of establishment with branching processes~\citep{athreya-ney-04}. This approach assumes that individuals survive, grow, and reproduce independently of one another and has provided fundamental insights into fixation of beneficial alleles~\citep{haldane-27}, the build up of biodiversity on islands~\citep{mcarthur-wilson-67}, the viability of endangered populations~\citep{soule-87}, and the evolution of disease emergence in novel host populations~\citep{antia-etal-03,park-etal-13}. However, even when populations are at low abundance, individuals may interact with one another (e.g. finding or competing for mates in sexually reproducing populations) and thereby violate the assumption of independence of these classical branching processes. When these interactions occur between different types of individuals, they lead to frequency-dependent feedbacks on the population dynamics.

To account for these frequency-dependent interactions within a founding population, \citet{Sch01} introduced a class of generalized urn models which were studied more extensively by \citet{BenSchTar04}. These models consider an urn containing a finite number of balls (the population) of different colors (the different genotypes or phenotypes). At each stage, balls of possibly different colors can be added or removed from the urn, modeling deaths, births, or changes of state due to interactions between individuals.

Two key questions about these Markov processes are a) when is there a positive probability that the population never goes extinct (i.e. the population establishes)? and b) on the event of non-extinction, what can be said about the long-term frequency dynamics? To address these questions,  \cite{Sch01} introduced a mean limit ordinary differential equation (ODE) on the simplex (corresponding to all possible population frequencies) associated with the urn models. Using these mean field ODEs \cite{BenSchTar04} proved: (i) if the population is expected to grow uniformly in the neighborhood of an attractor of this mean limit ODE, then with positive probability the population never goes extinct and the frequencies of the population converge to this attractor (see Theorem \ref{th:BST} next section); (ii) conversely if the population  is expected to decrease uniformly in the neighborhood of a given set, then convergence toward this set occurs with probability zero (see Theorem \ref{th:BST2} next section). As the expected growth rate of the population typically varies along non-equilibrium attractors, these two results, however, are most useful for equilibrium attractors of the mean limit ODE.

Here, we study the case where the underlying mean limit ODE admits non-equilibrium attractors with non-constant growth rates. As we show in the applications section, this case arises quite naturally in stochastic models for evolutionary games and population genetics. We extend the results of \cite{BenSchTar04} to a more general framework (see Theorems \ref{thm:attractor} and \ref{thm:attractor2} in Section \ref{sc:mr}). Most notably, we replace the assumption of uniform positive (respectively, negative) population growth near the attractor  with the assumption that the temporal average of the population growth rate is positive (respectively, negative) for initial conditions near the attractor (see assumptions (\ref{eq:PAG}) and (\ref{eq:NAG})).

The remainder of this paper is organized as follows. In next section, we define the class of generalized urn models, and recall the main results of \cite{BenSchTar04}. We also discuss the stochastic approximation methods, and briefly explain how they were used to derive these results. In section 3, we state and prove our main results: convergence with positive probability toward an attractor with average positive growth and non-convergence to an invariant set with negative average growth. 
Section 4 is devoted to applications to evolutionary games and population genetics. The proofs of some technical estimates are given in the Appendix.

\section{Generalized Urn Models} \label{sec:GUM}

In Section \ref{sc:model} we give the definitions of a class of generalized urn models  introduced by \citet{Sch01}.  The approach used to study these models is the so-called ODE method, which relates the asymptotic behavior of a stochastic difference equation to an ODE. This method is described in Section \ref{sc:ode}.

\subsection{The Urn Models} \label{sc:model}
We consider a finite population  consisting of individuals that are one of  $k$ types. Therefore, the state space for this Markov chain is the non-negative cone
    $$
    \Z^k_+=\{z=(z^1,\dots,z^k)\in \Z^k: z^i\ge 0\mbox{ for
     all } i \}
    $$
where $\Z^k$ is the space of $k$-tuples of integers. Given a vector
$w=(w^1,\dots,w^k)\in \Z^k$ define
    $$
    |w|=|w^1|+\dots+|w^k|\mbox{
    and }\alpha(w)=w^1+\dots+w^k.
     $$
We let $\|\cdot \|$ denote the  Euclidean norm on
$\R^k$.

Let $z(n)=(z^1(n),\dots,z^k(n))$ be a homogeneous Markov chain with
state space $\Z^k_+$. In our context, $z^i(n)$ corresponds to the number of  individuals of type $i$ at the $n$-th update.  Associated
with $z(n)$ is the random process $x(n)$ defined by
    $$
    x(n)=
    \left\{
    \begin{array}{ll}
    \frac{z(n)}{|z(n)|}&\mbox{if }z(n)\neq 0 \\
    &\\
    0&\mbox{if }z(n)=0
    \end{array}
    \right.
    $$
which is the distribution of types in the population at the $n$-th update. When there are  no individuals in the population at the $n$-th update, we arbitrarily set $x(n)$ to zero
which we view as the ``null'' distribution.  Provided, $z(n)$ is non-zero, the population distribution $x(n)$ lies on the  the unit  simplex
    $$
    S_k=\left\{ x=(x^1,\dots,x^k)\in \R^k: x^i\ge 0,\sum_{i=1}^k x^i=1\right\}.
    $$

We let $\Pi : \Z_+^k \times \Z_+^k \mapsto [0,1]$ denote the
transition kernel of the Markov chain $z(n)$. Specifically,
\[\Pi(z,z')=\P[z(n+1)=z'|z(n)=z].\]  Our standing
assumptions on the Markov chains $z(n)$ are as follows:
\begin{list}{}{}
\item (A1) At each update, there is a maximal number of individuals that can be
added or removed. In other words, there exists a positive integer
$m$ such that $|z(n+1)-z(n)|\le m$ for all $n$.
\item (A2) There exist Lipschitz maps $$\{p_w : S_k \to [0,1] \; :w \in \Z^k, \, |w| \leq m\}$$ and a real number $a>0$ such that $$ |p_w(z/|z|) -
\Pi(z,z+w)| \leq a/|z|$$ for all non-zero $z\in \Z^k_+$ and $w\in
\Z^k$ with $|w|\le m$.
\end{list}
Assumption (A2) implies that, as the population gets large, the transition probabilities tend to only depend  on the frequency vector $z/|z|$.

\subsection{Mean Limit ODEs} \label{sc:ode}

The following lemma which was proved by \cite{BenSchTar04} expresses the random process $(x(n))_n$ as a stochastic approximation algorithm. 

\begin{lemma}\label{lemma:gu}
Let $z(n)$ be a Markov chain on $\Z^k_+$ satisfying assumptions
(A1) and (A2) with mean limit transition probabilities $p_w:S_k\to
[0,1] $. Let ${\FF}_n$ denote the $\sigma$-field generated by $\{z(0), z(1), \ldots, z(n)\}.$ There exists sequences of random variables $\{U_n\}$ and $\{b_n\}$ adapted to ${\FF}_n$, and a real number $K>0$ such that
\begin{list}{}{}
\item (i) if $z(n)\neq 0$, then
    \begin{equation}\label{eq:sa}
    x(n+1)-x(n)=\frac{1}{|z(n)|}\left(\sum_{w\in \Z^k} p_w(x(n))(w-x(n) \alpha(w))+U_{n+1}+b_{n+1}\right).
    \end{equation}
\item (ii) $\E[U_{n+1}|z(n)]=0$.
\item (iii) The random variables $\| U_n\|$ and $\E[\|U_{n+1}\|^2|{\FF}_n]$ are uniformly bounded.
\item (iv) $\|b_{n+1}\|\le \frac{K}{\max\{1,|z(n)|\}}$.
\end{list}
\end{lemma}

The recurrence relationship (\ref{eq:sa}) can be viewed as a
``noisy'' Cauchy-Euler approximation scheme with step size
$1/|z(n)|$ for solving the ordinary differential equation
    \begin{equation}\label{eq:ave}
    \frac{dx}{dt}= g(x) := \sum_{w\in \Z^k} p_w(x)(w-x\alpha(w))
    \end{equation}

The limiting behavior of the $x(n)$ is therefore related to the solutions  of (\ref{eq:ave}). Indeed, when the number of individuals in the population grows without bound, the step size decreases to zero and it seems reasonable that there is a strong relationship between the limiting behavior of the mean limit ODE and the population distribution  $x(n)$. To make the relationship between the stochastic process $x(n)$ and the mean limit ODE more transparent, it is useful to define a continuous time version of $x(n)$ where time is scaled in an appropriate manner. Since the
number of events (updates) that  occur in a given time interval is
likely to be proportional to the size of the population, we define
{\em the time $\tau(n)$ that has elapsed by update $n$} as
    \begin{eqnarray*}
    \tau(0)&=&0\\
    \tau(n+1)&=&\left\{
    \begin{array}{ll}
    \tau(n)+\frac{1}{|z(n)|}&\mbox{if }z(n)\neq 0\\
    \tau(n) +1&\mbox{if }z(n)=0.
    \end{array}
    \right.
    \end{eqnarray*}
The continuous time version of $x(n)$ is given by
    \begin{equation}\label{eq:X}
    X(t)=x(n) \qquad \mbox{for }\tau(n)\le t<\tau(n+1).
    \end{equation}

To relate the limiting behavior of the solutions of (\ref{eq:ave}) to the limiting behavior of $X(t)$, we need a few more definitions to state a key result of \citet{Sch01}. Let $x.t$ denote the solution of \eqref{eq:ave} with initial condition $x$, at time $t$. A set $C$ is called {\em invariant} for \eqref{eq:ave} provided that $C.t = C$ (where $C.t = \left\{x.t: x \in C \right\}$)  for all $t\in \R$. A compact invariant set $\A\subset S_k$ is  an {\em attractor} if there is an open neighborhood
$U\subset S_k$ of $\A$ such that
    $$
    \cap_{t>0}\overline{\cup_{s\ge t} U.s}=\A.
    $$
The basin of attraction $B(\A)$ of $\A$ is the set of points $x\in
S_k$ satisfying $\inf_{y\in \A}\|x.t - y\|\to 0$ as
$t\to\infty$.
Finally, a compact invariant set $C$ is {\em
internally chain transitive} provided that \eqref{eq:ave} restricted to $C$ admits no proper attractor. 

Given a function $X:\R_+\to \R^k$ or a sequence $\{x(n)\}_{n\ge 0}$ in $\R^k$, we define the {\em limit sets},
$L(X(t))$ and $L(x(n))$, of $X(t)$ and $x(n)$ as follows. $L(X(t))$ is the set of
$p\in \R^k$ such that $\lim_{k\to\infty}X(t_k)=p$ for some
subsequence $\{t_k\}_{k\ge 0}$ with $\lim_{k\to\infty}t_k=\infty$.
$L(x(n))$ is the set of $p\in \R^k$ such that
$\lim_{k\to\infty}x(n_k)=p$ for some subsequence $\{n_k\}_{k\ge
0}$ with $\lim_{k\to\infty}n_k=\infty$.

Using the methods of \citet{Ben96, Ben99}, the following result of \cite{Sch01} demonstrates the relationship between the asymptotic behavior of $X(t)$ and $x.t$ on the event that the population is growing sufficiently rapidly. 

\begin{theorem}[Schreiber, 2001]\label{thm:apt}
Let $z(n)$ be a Markov process satisfying the assumptions of Lemma \ref{lemma:gu}. Then, on the event
\[\left\{ \sum_n \frac{1}{|z(n)|^{1+a}} < \infty, \; \mbox{ for some } a >0 \right\}, \]

\begin{enumerate}
\item the interpolated process $X(t)$ is almost surely an {\em asymptotic
pseudotrajectory} for the flow  of the mean limit ODE (\ref{eq:ave}). 
In other words, $X(t)$ almost surely satisfies
    \begin{equation} \label{eq:apt}
    \lim_{t\to \infty} \sup_{0\le h\le T}\|X(t).h - X(t+h)\|=0
    \end{equation}
for any $T>0$.
\item the limit set $L(X(t))$ of $X(t)$ is almost surely an internally chain transitive set for the mean limit ODE.
\end{enumerate}
\end{theorem}

The first assertion of the theorem roughly states that $X(t)$
tracks the solutions of the mean limit ODE (\ref{eq:ave}), with increasing accuracy far into the future. The second assertion of the theorem states that
the only candidates for limit sets of the population distribution $x(n)$
are connected compact internally chain recurrent sets for the mean limit ODE. With regards to attractors, we have the following useful property of internally chain recurrent sets (see, e.g., \citet[Cor. 5.4]{Ben99}) 
\begin{remarque} \label{rq:basinICT}
If an internally chain transitive set $C$ meets the basin of attraction of a given attractor $\A$ then it is contained in $\A$.  
\end{remarque}

Consider an attractor $\A$ for the mean limit ODE. Now suppose that there is a neighborhood $U\subset B(\A)$ of $\A$ such that whenever the population is large and its distribution $x(n)$ lies in $U$, the population is expected to grow. Specifically, $f(x(n))>0$ whenever $x(n)\in U$ where
\[
f(x)=\sum_w p_w(x) \alpha(w)
\]
is the limiting expected change in the population size. Under such circumstances, we would expect the population to increase in size and, consequently, the frequencies to follow the solutions of the mean limit ODE more closely. As $x(n)$ lies in the basin of attraction of $\A$, $x(n)$ would tend to remain near $\A$ and the population would be expected to increase further. One would expect that this positive feedback loop would result in the population growing with positive probability and its distribution converging to the attractor $\A$. Indeed, as the next theorem shows, this argument holds when  the attractor $\A$ is \emph{attainable}, which basically means that the random process can reach the basin of attraction of $\A$, at any time. More specifically, we define the set of {\em attainable points}, $Att_{\infty}(X)$, as
the set of points $x\in S_k$ such that, for all $M \in \mathbb{N}$ and
every open neighborhood $U$ of $x$
\[\mathbb{P}[|z(n)|\geq M \mbox{ and } x(n)\in U \mbox{ for some } n] > 0.\]

\begin{theorem}[Bena\"im, Schreiber and  Tarres, 2004] \label{th:BST} Let $z(n)$ be a generalized urn process verifying $(A1)$ and $(A2)$. Assume that 
\begin{equation}  \label{eq:unif_growth}
\lambda = \inf_{x \in \A} f(x) >0.
\end{equation}
If $B(\A) \cap Att_{\infty}(X) \neq \emptyset$, then
\[\mathbb{P} \left[ \liminf_n \frac{|z(n)|}{n} \geq \lambda  \, \mbox{ and } \;  L(x(n)) \subset  \A \right] >0\]
\end{theorem}

The next theorem provides a partial converse to Theorem~\ref{th:BST}. Roughly, it states that if there is a compact set $K\subset S_k$ near which the population is expected to decrease every update, then the population distribution $x(n)$ can not converge to $K$. 

\begin{theorem}[Bena\"im, Schreiber and  Tarres, 2004] \label{th:BST2} Let $z(n)$ be a generalized urn process verifying $(A1)$ and $(A2)$ and $K \subset S_k$ be any compact set. Assume that 
\begin{equation}  \label{eq:unif_neg_growth}
 \sup_{x \in K} f(x) <0.
 \end{equation}
Then there exists $M>0$ such that 
\[\mathbb{P} \left(  |z(n)| >M \; \mbox{ for n large enough} \, \mbox{ and } \; L(x(n)) \subset  K \right) =0\]
\end{theorem}

\section{Main results} \label{sc:mr}

While the assumption of uniform growth is not restrictive for equilibrium attractors of the mean limit ODE, the long-term behavior of these mean limit ODEs may be governed by non-equilibrium behavior, such as limit cycles, quasi-periodic motions, or chaotic attractors. For these types of attractors, the growth rate $f(x)$ of a population typically varies along points of the attractor and, consequently, the uniform growth assumptions (\ref{eq:unif_growth}) and (\ref{eq:unif_neg_growth}) are too restrictive. Throughout the section, ${\mathcal A}$ denotes an attractor for the mean limit ODE, with basin of attraction $B({\mathcal A})$. Rather than a uniform growth assumption, we assume that the long-term temporal average of the growth rate $f$ is positive along orbits of the mean limit ODE. Specifically, 
\begin{equation}  \label{eq:PAG}
 \liminf_{t \rightarrow + \infty} \frac{1}{t} \int_{0}^t f(x.s) ds >0 \mbox{ for all } x \in B(\mathcal{A}). 
\end{equation}

Our main result is a generalization of Theorem \ref{th:BST}, under this less restrictive positive growth assumption (\ref{eq:PAG}). This result states, roughly, that if the basin of attraction is attainable and the temporal averages of the growth rate are positive along the attractor, then the population grows without bound and its distribution converges to the attractor with positive probability. 

\begin{theorem}[Positive average growth and convergence] \label{thm:attractor} Let $z(n)$ be a generalized urn
process satisfying (A1) and (A2). Assume that (\ref{eq:PAG}) holds
and 
$$B({\mathcal A})\cap Att_{\infty}(X)\neq
\emptyset,
$$ 
then 
\[\P \left[\sum_{n} \frac{1}{|z(n)|^{1+ \delta}} < + \infty \, \,  \forall \delta>0, \; \mbox{ and } \,  L(x(n)) \subset \mathcal{A} \right]>0.\]
\end{theorem}

Unlike Theorem~\ref{th:BST}, Theorem~\ref{thm:attractor} can no longer guarantee  linear population growth with positive probability. However, this is not surprising, as the population growth rate can be negative as well as positive despite having a positive temporal average. 

We have a partial converse result to Theorem~\ref{thm:attractor}. We say that there is an average negative growth rate in an invariant compact set $K$ if
\begin{equation}  \label{eq:NAG}
\limsup_{t \rightarrow + \infty} \frac{1}{t} \int_{0}^t f(x.s) ds <0 \mbox{ for all }x\in K.
\end{equation}

\begin{theorem}[Negative average growth and non-convergence] \label{thm:attractor2} Let $z(n)$ be a generalized urn
process satisfying (A1) and (A2). Assume that (\ref{eq:NAG}) holds for a compact invariant set $K$. Then there exists $M>0$ such that 
\[\mathbb{P} \left[ |z(n)| \geq M \; \mbox{ for n large enough} \, \mbox{ and } \;  L(x(n)) \subset K  \right]=0.\]
\end{theorem}

The proofs of Theorems~\ref{thm:attractor} and \ref{thm:attractor2} are given in Sections 3.2 and 3.3, respectively. Several key technical estimates required for these proofs are described in 3.1 and proven in the Appendices. 

\subsection{Key estimates\label{sec:technical}} 

To state the key estimates for the proofs of our main results, 
call $L_f$ the Lipschitz constant of $f$ and $\|f\|_{\infty} :=\sup_x f(x)$.  The map $s \mapsto f(x.s)$ is Lipschitz, uniformly in $x \in K$. Consequently, there exists a constant $L'>0$ such that
\[\|f(x.s) - f(x.s')\| \leq L'|s-s'| \, \mbox{ for all } x \in K.\]
Recall that $g(x)=\sum_w p_w(x)(w-x\alpha(w))$. Let $L_g$ be the Lipschitz constant for $g$ and $\|g\|_{\infty} =\sup\|g(x)\|$. Define $L= \max \{L_f,L_g,L',\|g\|_{\infty},\|f\|_{\infty}\}$.

Now assume that $(x(n))_n$ is a stochastic approximation process which satisfies equation \eqref{eq:sa}.  Define  $\|U\| = \sup_n \|U_n\|$. To simplify the presentation of the proof, we assume that $\|b_n\|=0$. The proof without this assumption is notationally more cumbersome but follows in a nearly identical fashion.

Define $m(t) = \inf \{k : \; \, \tau(k) \geq t\}$. Let  $T_0>0$.  We first observe that the population size $|z(r)|$ remains bounded on time intervals of order $T_0$. Given $r_0\in \N$, define $r_{k+1}= m(\tau(r_k)+T_0)$ for all $k\ge 1$. 
Provided $|z(r_k)|$ is large enough, \citet[Lemma 2]{BenSchTar04} states that
\begin{enumerate} 
\item[$(i)$] $B|z(r_k)|\geq|z(r)|\geq B^{-1}|z(r_k)|$ for all $r\in [r_k,r_{k+1}]$, and
\item[$(ii)$] $T_0 B^{-1} |z(r_k)|\leq r_{k+1}-r_k\leq T_0 B |z(r_k)|$
\end{enumerate}
where $B=3 e^{mT_0}$. Using Gronwall's inequality, \citet[Proposition 4.1]{Ben99} proved the following estimate 
\begin{equation} \label{eq:erg}
\sup_{r \in [r_k,r_{k+1}]}\|x(r_k).(\tau(r) - \tau(r_k)) - x(r)\|\leq C \left(
\Gamma_1(k,T_0)+\Gamma_2(k,T_0)\right)
\end{equation}
where $C$ is a positive constant, which depends on $T_0$ and $L$, 
$$
\Gamma_1(k,T_0)=\sup_{r \in [r_k,r_{k+1} -1]} \left|
\sum_{i=r_k}^{r} \frac{U(i+1)}{|z(i)|}\right|
$$
and \footnote{in the more general case where $\|b\| \neq 0$, there is an additional term in $\Gamma_2$, namely $\sup_{r \in [r_k,r_{k+1} -1]} \left|
\sum_{r_k}^{r} \frac{b(i+1)}{|z(i)|} \right| 
$}
$$
\Gamma_2(k,T_0)=\frac{2\|g\|_{\infty} }{ \inf_{r \in [r_k,r_{k+1}] }|z(r)|}$$

The next Lemma refines the statement of the asymptotic pseudotrajectory property (\ref{eq:apt}) for the discrete time process $x(n)$. The proof is given in the Appendix.

\begin{lemma} \label{lm:diff} Let $T_0$ and $\delta$ be positive real numbers. Then we have
\begin{equation} \label{eq:upbound}
\mathbb{P} \left(\sup_{r \in [r_k,r_{k+1}]}\|x(r_k).(\tau(r) - \tau(r_k)) - x(r)\| > \delta \Big| \; \, \mathcal{F}_{r_k}\right) \leq  \frac{C_0(T_0)}{|z(r_k)| \delta^2},
\end{equation}
(where $C_0(T_0) := 4\|U\|^2BC^2 T_0$) on the event $\left\{|z(r_k)| \geq \frac{4BCL}{\delta} \right\}.$

\end{lemma}

The next two propositions roughly underestimate the likelihood that the average growth $\frac{1}{T_0} \int_0^{T_0} f(x(r).s) ds$ remains close of its stochastic counterpart on intervals of time of length $T_0$, provided the population size is initially large enough. Both proofs are given in the Appendix.

\begin{proposition} \label{pr:tech} Let $(x(n))_n$ be a stochastic process satisfying (\ref{eq:sa}). Given a point $y\in S_k$, a time $T_0>0$, and $\delta>0$, there exists a compact neighborhood $U$ of $y$, $C_1(T_0)>0$, and $M_1>0$ such that
\begin{equation}\label{eq:tech}
\mathbb{P} \left[\frac{1}{T_0} \left| \sum_{r=r_k}^{r_{k+1}-1} \frac{1}{|z(r)|} f(x(r))  -\int_0^{T_0} f(y.s)ds\right| >\delta\Big|\mathcal{F}_{r_k}\right] \leq \frac{C_1(T_0)}{|z(r_k)|\delta^2}\end{equation}
on the event   $V_k := \left\{x(r_k) \in U \right\} \cap \left\{|z(r_k)| \geq M_1\right\}$.
\end{proposition}

Using the estimate from this proposition, we can get the following result.

\begin{proposition} \label{pr:prop2} Let $(x(n))_n$ be a stochastic process satisfying (\ref{eq:sa}). Given a point $y\in S_k$, a time $T_0>0$, and $\delta>0$, there exists a compact neighborhood $U$ of $y$, $C_2(T_0)>0$, and $M_2>0$ such that
\[
\mathbb{P}\left[\frac{1}{T_0}\left|  \sum_{i=r_k}^{r_{k+1}-1} \frac{|z(i+1)|-|z(i)|}{|z(i)|}-\int_0^{T_0} f(x(r_k).s)ds\right| \ge \delta \Big| \mathcal{F}_{r_k} \right] \le \frac{C_2(T_0)}{|z(r_k)|\delta^2}
\]
on the event $W_k:= \left\{x(r_k) \in U \right\} \cap \left\{|z(r_k)| \geq M_2\right\}$.
\end{proposition}
\vspace{.2cm}

\subsection{Proof of Theorem \ref{thm:attractor}} 

Pick $p \in B(\mathcal{A}) \cap Att_{\infty}(X)$ and an open neighborhood $U$ of $A$, which contains $p$ and whose closure $K$ is compact and included in $B(\mathcal{A})$. Since $\mathcal{A}$ is an attractor, there exists a positive time $T'$ and $\delta>0$ such that, for any $T \geq T'$, 
\[X(t) \in U \mbox{ and } \;  \|X(t+T)-X(t).T\| < \delta \Rightarrow X(t+T) \in U.\]
Also, by assumption (\ref{eq:PAG}), there exists $a_1 >0$ and $T'' >0$ such that,  for any $T \geq T'' $ and $x \in K$,
\[\frac{1}{T} \int_0^{T} f(x.s) ds \geq a_1.\]
Define $T_0 = \max \{T',T'' \}$.

Let $M>0$ be fixed (we will need $M$ to be greater than some quantity which depends on $T_0$, $m$, $a_1$ and $\delta$ in a manner that will be specified later in the proof). By the attainability condition, there exists an index $r_0 \in \mathbb{N}$ such that
\[\mathbb{P} \left[x(r_0) \in U, \; \, |z(r_0)| \geq M \right] >0.\]
Consider the following events for all $k\ge 1$
$$
E_1(k) = \{|z(r_k)|\geq \zeta^{k-1} B^{-1} M\}; \; \; E_2(k) = \{x(r) \in U, \mbox{ for all } r\in[r_{k},r_{k+1}]\}$$
where $\zeta=1+a_1T_0/4$. Let $E(0)$ be the event $\{ |z(r_0)|\ge
M, x(r_0) \in U \}$.  For $k\geq 1$, define $E(k)=E(k-1)\cap
E_1(k)\cap E_2(k)$. We will show that there exists a constant
$F>0$ such that 
\begin{equation}\label{eq:key}
\mathbb{P}[E(k+1) | E(k)] \geq 1 - F/M\zeta^k \mbox{ for all } k\ge 0.
\end{equation} 

The proof of estimate \eqref{eq:key} relies on two ingredients. First, Lemma~\ref{lemma2}  provides a lower bound to  the probability of being inside $U$ on $[r_{k+1},r_{k+2}]$ if $x(r_k)\in U$. Second, Lemma~\ref{lemma3} underestimates  the probability that the population grows sufficiently (namely by a multiplicative parameter $\zeta>1$) between times $r_k$ and $r_{k+1}$, provided it stayed in $U$ the entire  
 time. Unlike the work of \citet{BenSchTar04}, the main issue here is that  we do not have expected growth for every update of the population in the neighborhood of the attractor, so we need to make use of Proposition \ref{pr:prop2} and the estimates on the population size given above.

\begin{lemma}\label{lemma2} There exist  
$M'_0 >0$ and $D>0$  such that if $M\ge M_0'$, then 
$$P[E_2(k+1) | E(k)] \geq 1-\frac{D}{M\zeta^{k-1}}$$
for all $k\ge 0$.\end{lemma}

\noindent {\bfseries Proof.}
Statement \emph{(i)} in Section~\ref{sec:technical} implies that $|z(r)|\geq B^{-2}|z(r_{k})|$ for $r\in[r_{k}, r_{k+2}]$. Furthermore, the definition of $r_k$ implies that 
 $T_0\leq \tau(r_{k+1})-\tau(r_k)\leq 2 T_0$ which implies $4T_0\geq
\tau(r)-\tau(r_{k})\geq T_0$ for $r\in[r_{k+1}, r_{k+2}]$. Statements \emph{(i)} and \emph{(ii)} in Section~\ref{sec:technical} also imply that $r_{k+2}-r_{k}\le 2 B^2T_0|z(r_{k})|$ and  $r_{k+2} \leq m(\tau(r_k) + 3T_0)$.

On the event $E_1(k)$, $|z(r)|\ge |z(r_k)|B^{-2}\ge\zeta^{k-1}B^{-3}M$ for all $r\in[r_{k},r_{k+2}]$. Therefore, we can apply Lemma \ref{lm:diff} with $3T_0$ and $|z(r_k)| \geq \zeta^{k-1}B^{-3}M$ which gives, choosing $M'_0$ large enough and $M$ greater than $M'_0$,
\begin{eqnarray*}
\mathbb{P} \left[\sup_{r\in[r_{k+1}, r_{k+2}]}
\left\| x(r_{k}).(\tau(r)-\tau(r_k)) - x(r) \right\| > \delta \mid \; \, \mathcal{F}_{r_k} \right] \leq \frac{B^3C_0(3T_0)}{\delta^2 M\zeta^{k-1}}.
\end{eqnarray*}
Since $\tau(r_{k+1})-\tau(r_{k})\ge T_0$ and $x(r_k)\in U$, our
choice of $T_0$ and $\delta$ implies that $x(r)\in U$ for all $r\in
[r_{k+1},r_{k+2}]$, on the event \[\left\{\sup_{r\in[r_{k+1}, r_{k+2}]}
\left\|x(r_{k}).(\tau(r) - \tau(r_k)) - x(r) \right\|\le \delta\right\}.\] Defining $D = B^3C_0(3T_0)/ \delta^2$, we therefore have
\[\mathbb{P} [E_2(k+1)^c \mid \; E(k) ] \leq \frac{D}{M\zeta^{k-1}} \]
$ \blacksquare$

\begin{lemma}\label{lemma3}
There exist $M''_0 \geq M'_0$ and $D'>0$ such that if $M \geq M_0''$, then
$$P[E_1(k+1) \mid \;  E(k)]\geq 1-\frac{D'}{M\zeta^{k-1}}$$
for all $k\ge 1$.
\end{lemma}

\noindent {\bfseries Proof.} On the event $E_1(k+1)^c \cap E(k)$, we have 
\[
n_k := |z(r_{k+1})| - |z(r_k)| \le \zeta^k B^{-1}M -\zeta^{k-1} B^{-1}M = B^{-1}M \zeta^{k-1}(a_1T_0/4).\] Under the constraint that one cannot add more balls than this fixed quantity  between times $r_k$ and $r_{k+1}$, we prove the following inequality in the Appendix
\begin{equation} \label{claim}  
\sum_{i=r_k}^{r_{k+1}-1} \frac{|z(i+1)|-|z(i)|}{|z(i)|} \leq   \frac{n_k}{|z(r_k)|} + 
 \frac{m}{2}(r_{k+1} - r_k) \left(\frac{1}{B^{-1}|z(r_k)|} - \frac{1}{B^{-1}|z(r_k)|+m}\right).
\end{equation}
Consequently, since $r_{k+1}-r_k \leq T_0B|z(r_k)|$,
\begin{eqnarray*}
\sum_{i=r_k}^{r_{k+1}-1} \frac{|z(i+1)|-|z(i)|}{|z(i)|} &\leq&  B^{-1}M \zeta^{k-1}(a_1T_0/4) \frac{1}{|z(r_k)|} + \frac{T_0 B |z(r_k)|m^2}{2(B^{-2}|z(r_k)|^2)}\\
&\leq& \frac{a_1T_0}{4} + \frac{T_0 B^3 m^2}{2M} .
\end{eqnarray*}
Hence, choosing $M$ greater than $M_0'' := \max \{M_0',2T_0B^3m^2/a_1\}$, 
\[\frac{1}{T_0}\sum_{i=r_k}^{r_{k+1}-1} \frac{|z(i+1)|-|z(i)|}{|z(i)|} < \frac{a_1}{2}.\]

Recall that, by definition of $T_0$, 
\[
\frac{1}{T_0}\int_0^{T_0} f(x(r_k).s)ds \geq a_1\]
Hence, on the event $E_1(k+1)^c \cap E(k)$, we have
\[\frac{1}{T_0}\left|  \sum_{i=r_k}^{r_{k+1}-1} \frac{|z(i+1)|-|z(i)|}{|z(i)|}-\int_0^{T_0} f(x(r_k).s)ds\right| \geq \frac{a_1}{2}.\]

Finally, by Proposition \ref{pr:prop2} and provided $M''_0$ is large enough, we have
\[\mathbb{P} \left[E_1(k+1)^c \mid E(k) \right] \leq \frac{4C_2(T_0)}{a_1^2|z(r_k)|}\]
The proof is complete, taking $D' = \frac{4C_2(T_0)B}{a_1^2}$. $\blacksquare$

Now choose $M$ larger than $M_0''$. By Lemmas \ref{lemma2} and \ref{lemma3}, and denoting $F = D +D' >0$,  we have 
$$\P[E(k+1)|E(k)]\ge 1-\frac{F}{M\zeta^{k-1}}$$
for all $k\ge 0$. Since the sequence of event $\{E(k)\}_k$ is decreasing, 
it follows that
\begin{eqnarray*}
\P\left[\lim_{k\to\infty} E(k)\right]&\ge&\P[E(0)]
\left(1-\sum_{k=1}^{\infty}\frac{F}{M\zeta^{k-1}}\right)\\
&\ge& \P[E(0)]\left(1-\frac{\zeta F}{M(\zeta-1)}\right).
\end{eqnarray*}

On the event $\lim_{k\to\infty} E(k)$, we have for any $\delta>0$
\begin{eqnarray*}
\sum_{i=1}^{+ \infty} \frac{1}{|z(i)|^{1+\delta}} &=& \sum_{i=1}^{r_0} \frac{1}{|z(i)|^{1+\delta}} + \sum_{k =1}^\infty \sum_{i=r_k+1}^{r_{k+1}} \frac{1}{|z(i)|^{1+ \delta}}\\
&\le &\sum_{i=1}^{r_0} \frac{1}{|z(i)|^{1+\delta}} + \sum_{k =1}^\infty \left(\max_{r_k+1\le i \le r_{k+1}}\frac{1}{|z(i)|^{\delta}} \times \sum_{i=r_k+1}^{r_{k+1}} \frac{1}{|z(i)|} \right)\\
&\le & \sum_{i=1}^{r_0} \frac{1}{|z(i)|^{1+\delta}} + T_0 \sum_{k =1}^\infty \max_{r_k+1\le i \le r_{k+1} }\frac{1}{|z(i)|^{\delta}}\\
&\le & \sum_{i=1}^{r_0} \frac{1}{|z(i)|^{1+\delta}} + T_0 \sum_{k =1}^\infty  \frac{1}{(B^{-1} |z(r_k)|)^{\delta}}\\
&\le & \sum_{i=1}^{r_0} \frac{1}{|z(i)|^{1+\delta}} + T_0  \sum_{k =1}^\infty  \frac{1}{( \zeta^{k-1} B^{-2}M)^{\delta}}\\
&=& \sum_{i=1}^{r_0} \frac{1}{|z(i)|^{1+\delta}} + T_0 B^{2\delta} M^{-\delta} \frac{\zeta^{\delta}}{\zeta^{\delta}-1} <\infty.\\
\end{eqnarray*}

Hence the definition of $E_2(k)$ implies that $$\P[{\mathcal C}]\ge \P[E(0)]\left(1-\frac{\zeta F}{M(\zeta-1)}\right) >0,$$ 
where 
$${\mathcal C}=\left\{\sum_{i=1}^{+\infty} \frac{1}{|z(i)|^{1+\delta}} < \infty, \; \forall \delta>0 \; \mbox{ and } \; \;  x(n) \in U \, \forall n\right\}.$$
On the event ${\mathcal C}$, Theorem~\ref{thm:apt} implies that  $L(\{x(n)\})$ is a compact
internally chain recurrent set for the mean limit ODE. Since
$L(\{x(n)\})\subset K \subset  B({\mathcal A})$ on the event
${\mathcal C}$, Remark \ref{rq:basinICT} and Theorem \ref{thm:apt} imply that
$L(\{x(n)\})\subset{\mathcal A}$. $\; \; \square$
\vspace{.2cm}

\subsection{Proof of Theorem \ref{thm:attractor2}} By assumption (\ref{eq:NAG}), there exists an open neighbourhood $U$ of $K$, $T_0 >0$ and $a_1 >0$ such that
\[
\frac{1}{T_0} \int_{0}^{T_0} f(x.s) ds < -a_1 \mbox{ for all }x \in U.
\]
Let $M_2$ given by Proposition~\ref{pr:prop2} and choose $M > \max\{ M_2, |z(0)|\}$. Define 
\[
T_\ell = \inf \left\{n \geq \ell: x(r_n) \notin U \,\mbox{ or } \,  |z(r_n)|<M \right\}.
\]
Then 
\begin{equation}\label{eq:Tl} \left\{L((x(n))_n) \subset K, \; \; |z(n)| \geq M \; \mbox{ for $n$ large enough} \right\} \subset \bigcup_{\ell \in \mathbb{N}^*} \{ T_\ell = + \infty\}.\end{equation} 
 For large enough $M$, we  prove that   $\mathbb{P}[T_\ell = +\infty] = 0$ for all $\ell$.

By Proposition \ref{pr:prop2},  there exists $\alpha >0$ which depends on $T_0$ and $a_1$ such that
\[
\mathbb{P}\left[\frac{1}{T_0}  \sum_{i=r_k}^{r_{k+1} -1} \frac{|z(i+1)|-|z(i)|}{|z(i)|} \ge -a_1/2  \right] \le \frac{\alpha}{|z(r_k)|}
\]
for any $k \geq \ell$, on the event $T_{\ell} > k$.  Let  $N(i) = |z(i+1)|-|z(i)|$. On the event $\sum_{i=r_k}^{r_{k+1}-1} \frac{N(i)}{|z(i)|}<0$, we claim that 
{\small
\begin{equation} \label{eq:9}
\sum_{i=r_k}^{r_{k+1}-1} \frac{N(i)}{|z(i)|} \geq \frac{|z(r_{k+1})| - |z(r_k)|}{B^{-1}|z(r_k)|+1}.
\end{equation}
}
To prove this inequality, let us first assume that by contradiction $|z(r_{k+1})| - |z(r_k)| \geq 0$. Then in the interval of time from $r_k$ to $r_{k+1}$, more balls were added than removed. Therefore, for every  ball removed when the number of balls is $|z(i)|$ for some $r_k \le i \le r_{k+1}$, there is  ball added when the number of balls is at most $|z(i)|-1$. Consequently $\sum_{i=r_k}^{r_{k+1}-1} N(i)/|z(i)| >0$, a contradiction. Hence $|z(r_{k+1})| - |z(r_k)| <0$: more balls are removed than added. For every ball that is added when the number of balls is $|z(i)|$ for some $r_k \le i \le r_{k+1}$, a ball is removed when the state is at least $|z(i)|+1$.  Now the remaining $|z(r_k)| - |z(r_{k+1})|$ balls were removed when the number of balls was at least $B^{-1} |z(r_k)|+1$. This proves (\ref{eq:9}). As a consequence 
\begin{eqnarray*}
|z(r_{k+1})| - |z(r_k)| &\leq& B^{-1}|z(r_k)| \sum_{i=r_k}^{r_{k+1}-1} \frac{N(i)}{|z(i)|}
\end{eqnarray*}
and
\[|z(r_{k+1})| - |z(r_k)| \leq B^{-1} |z(r_k)| \frac{-a_1T_0}{2}.\]
with probability greater than $1-\alpha/|z(r_k)|$ on the event $\{T_\ell >k\}$. Moreover, $|z(r_{k+1})| - |z(r_k)|$ can never be larger than $r_{k+1} - r_k < T_0 B |z(r_k)|$. Hence, 
{\small
\begin{eqnarray*}
\mathbb{E} \left[|z(r_{k+1})| - |z(r_k)| \mid T_\ell >k \right] &\leq&   \frac{-a_1T_0 |z(r_k)|}{2B} \left(1 - \frac{\alpha}{|z(r_k)|}\right) + T_0 B |z(r_k)| \frac{\alpha}{|z(r_k)|}\\
&\leq&  |z(r_k)| \frac{-a_1 T_0 }{2B} + \frac{\alpha a_1 T_0}{2B} + T_0 B \alpha\\
&\leq& -\epsilon |z(r_k)| \leq - \epsilon M.
\end{eqnarray*}
}for some $\epsilon >0$. The remainder of the proof is similar to the proof of \citet[Proposition 2]{BenSchTar04}: for $n > \ell$ we have 
\begin{eqnarray*}
0 \leq \mathbb{E}(|z(r_{n \wedge T_\ell})|)&=&  \sum_{k=l+1}^{n}  \mathbb{E}[|z(r_{k \wedge T_\ell})|-|z(r_{(k-1) \wedge T_\ell })|]+\mathbb{E}(|z(r_\ell)|)\\
&=& \sum_{k=l+1}^{n}  \mathbb{E}\left[|z(r_{k})|-|z(r_{k-1})| \mid T_\ell \geq k\right] \mathbb{P}[T_\ell \geq k) +\mathbb{E}(|z(r_\ell)|] \\
&\le & -\epsilon M \sum_{k=l+1}^n \mathbb{P}[T_\ell\ge k] +\mathbb{E}[|z(r_\ell)|]
\end{eqnarray*}
Taking the limit as $n\to\infty$, we get that $\sum_{k=\ell+1}^\infty
\mathbb{P}[T_\ell\ge k] \le \mathbb{E}[|z(r_\ell)|]/\epsilon M.$ The Borel-Cantelli lemma implies
that $\mathbb{P}[T_\ell=\infty]=0$. Relation (\ref{eq:Tl}) implies that
\[\mathbb{P}[\{L(x(n))\subset K\}\cap \{|z(n)| \geq M \; \, \mbox{for $n$ large enough} \} ]=0. \; \; \blacksquare \]

\section{Applications}

To illustrate the applicability of our results, we examine stochastic analogs of  two well-know evolutionary dynamics: replicator dynamics from evolutionary game theory~\citep{schuster-sigmund-83,hofbauer-sigmund-98} and selection-mutation dynamics from population genetics~\citep{hofbauer-85,bulmer-91,hofbauer-sigmund-98}.

\subsection{Replicator processes}

Consider a population consisting of individuals playing $k$ different strategies. In the absence of interactions, each individual produces offspring at rate $b$ and dies at  rate $d$. Each individual initiates an encounter with another individual at rate $\nu$. Individual encounters are random. If an individual with strategy $i$ initiated an encounter with an individual with strategy $j$, then the individual with strategy $i$ either gives birth with probability $b_{ij}$, dies with probability $d_{ij}$ or is unaffected by the encounter with probability $u_{ij}=1-b_{ij}-d_{ij}$. Similarly, the individual with strategy $j$ gives births, dies, or is unaffected with probabilities $b_{ji}$, $d_{ji}$, and $u_{ji}$.  
While this description yields a continuous-time Markov chain $\tilde z(t)$, we focus on the embedded discrete-time Markov chain $z(n)$  corresponding to the state of $\tilde z(t)$ at the $n$-th update via a birth, death, or encounter between two individuals. We note that the probability of unbounded population growth (i.e. non-extinction) and convergence of the population distribution to a particular set in $S_{k}$ are equal for the embedded and continuous-time processes. Hence, for the questions we are interested in, there is no loss of information by restricting our attention to the discrete-time model.

Let $e_1=(1,0,0,\dots,0), e_2=(0,1,0,0,\dots,0), \dots, e_k=(0,\dots,0,0,1)$ be the standard basis of $\mathbb{R}^k$. At each update,  $z(n)$ can change either by one individual giving birth ($+e_i$) possibly following an encounter, one individual dying ($-e_i$) possibly following an encounter, or pairs of births or deaths  following an encounter between individuals  ($e_i-e_j$, $-e_i-e_j$ or $e_i+e_j$ for $i,j \in \{1,2,\dots,k\}$). Setting $\gamma=\frac{1}{b+d+\nu}$, the limiting transition functions $p_w$, as $|z|\to\infty$, are defined by  
\begin{eqnarray*}
p_{e_i}(x) &=& \gamma x_i\left(b+ 2\nu \sum_{j=1}^k x_j b_{ij} u_{ji}\right)\\
p_{e_i+e_j}(x) &=& 2\gamma \nu x_i x_j b_{ij}b_{ji} \mbox{ for } i \neq j\\
p_{2e_i}(x)&=&\gamma \nu x_i^2 b_{ii}^2  \\
p_{e_i-e_j}(x)&= & 2 \gamma \nu  x_i x_j b_{ij}d_{ji}\mbox{ for } i\neq j\\
\end{eqnarray*}
and analogously for $p_{-e_i}$, $p_{-e_i-e_j}$ and $p_{-2e_i}$.
The transition functions, $\Pi(z,z+w)$, for the non-limiting process $z(n)$ are given by $p_w(x=\frac{z}{|z|})$ with the exception that 
\[\Pi(z,z+2e_i) = \gamma \nu x_i \frac{(x_i |z| -1)}{|z|} b_{ii}^2\]
(and analogously for $-2e_i$) to ensure that these interactions are between \emph{different} individuals playing the same strategy $i$.

Let $B=(b_{ij})_{i,j}$ and $D=(d_{ij})_{i,j}$ be the birth and death matrices. The mean limit ODE for the process is given by 
\begin{equation}
\frac{dx}{dt} = \sum_{w} p_w(x)(w-x\alpha(w)) = 2\nu \gamma \, x \circ \left( (B-D)x-x^T (B-D)x\right)
\end{equation}
where $\circ$ denotes the Hadamard product. Setting $A=2\nu\gamma(B-D)$ yields the standard form of the replicator equations \citep{hofbauer-sigmund-98}:
\begin{equation}\label{eq:replicator}
\frac{dx}{dt} = x\circ( Ax -x^T Ax)
\end{equation}
Two key results for these equations described by \citet{hofbauer-sigmund-98} are the following.

\begin{theorem}\label{thm:hs1} Let $A$ be a $k\times k$ matrix. 
If there exist $p_1,\dots,p_k>0$ such that 
\[
\sum_{i=1}^k p_i \left((Ax)_i -x^TAx \right)>0
\]
for all equilibria of (\ref{eq:replicator}) lying on the boundary $\partial S_k$ of the simplex, then there exists an attractor $\mathcal A$ for (\ref{eq:replicator}) whose basin of attraction is $S_k \setminus \partial S_k$.
\end{theorem}

\begin{theorem}\label{thm:hs2}  Let $A$ be a $k\times k$ matrix. If there exists a compact invariant set $K\subset S_k \setminus \partial S_k$ for (\ref{eq:replicator}), then there exists a unique positive equilibrium $\hat x \in S_k \setminus \partial S_k$ such that 
\[
\lim_{t\to\infty} \frac{1}{t}\int_0^t (x.s) \,ds =\hat x
\]
and 
\[
\lim_{t\to \infty}\frac{1}{t}\int_0^t (x.s)^T A (x.s)\,ds = \hat x^T A \hat x
\]
for any solution of (\ref{eq:replicator}) with initial condition in $K$.
\end{theorem}

These two Theorems in conjunction with our main results Theorems~\ref{thm:attractor} and \ref{thm:attractor2}, imply the following result. We recall that $F$ is a \emph{face} of $S_k$ if there exists a set $I\subset \{1,\dots, k\}$  such that $F=\{x\in S_k: x_i =0$ \mbox{ for }\;  $i\in I\}$.

\begin{theorem}  Let $z(n)$ be a replicator process with $b>0$ and $d>0$ and $\prod_i z_i(0)>0$. Define $A=2\nu\gamma(B-D)$ as above. Let $F$ be a face of the simplex $S_k$.  If there exist $p_1,\dots,p_k>0$ such that 
\begin{equation}\label{eq:one}
\sum_{i\notin I} p_i \left((Ax)_i -x^TAx \right)>0
\end{equation}
for all equilibria $x\in \partial F$ and 
\begin{equation}\label{eq:two}
\frac{b-d}{b+d+\nu} +   \hat x^T A \hat x>0
\end{equation}
for the interior equilibrium $\hat x\in F\setminus \partial F$, then there exists a compact set $\mathcal A$ in $F\setminus \partial F$ such that
\[
\mathbb{P}\left[ \sum_n \frac{1}{|z(n)|^{1+\delta}}<\infty \mbox{ for all } \delta>0 \mbox{ and } L(\{x(n)\}_n)\subset \mathcal{A}\right]>0.
\]   
Alternatively, if $K$ is a compact, invariant subset of $F\setminus \partial F$ and 
\begin{equation}\label{eq:three}
\frac{b-d}{b+d+\nu} +  \hat x^T A \hat x<0
\end{equation}
for the interior equilibrium $\hat x\in F\setminus \partial F$, then there exists $M>0$ such that
\[
\mathbb{P}\left[ |z(n)| \ge M \mbox{ for $n$ sufficiently large and } L(\{x(n)\}_n)\subset \mathcal{K}\right]=0.
\]   
\end{theorem}

\noindent \textbf{Proof:} To prove the first assertion, assume that $F$ is a face of the simplex $S_k$ such that \eqref{eq:one} holds. Since $\prod_i z_i(0)>0$ and $d>0$, there is a positive probability that $z(n)$ lies in the interior of the face $F$ at some update $n\ge 1$. Hence, without loss of generality and by shifting time, we assume that $z(0)$ lies in the interior of $F$.  Theorem~\ref{thm:hs1} implies that there is an attractor $\A$ in the interior of $F$ for the replicator dynamics restricted to $F$ and $B(\A)$ equals the interior of $F$. The assumptions that $b>0$ and $d>0$ imply that $B(\A)$ is attainable. Assumption~\eqref{eq:two} and Theorem~\ref{thm:hs2} imply there is average positive growth rate in $B(\A)$ i.e. \eqref{eq:PAG} holds.  Applying Theorem~\ref{thm:attractor} completes the proof of the first assertion.   

To prove the second assertion, assumption~\eqref{eq:three} and Theorem~\ref{thm:hs2} imply there is average negative growth rate at $K$ i.e. \eqref{eq:NAG} holds. Applying Theorem~\ref{thm:attractor2} completes the proof of the second assertion.   \qed
\vspace{.2cm}

\begin{figure}
\includegraphics[width=\textwidth]{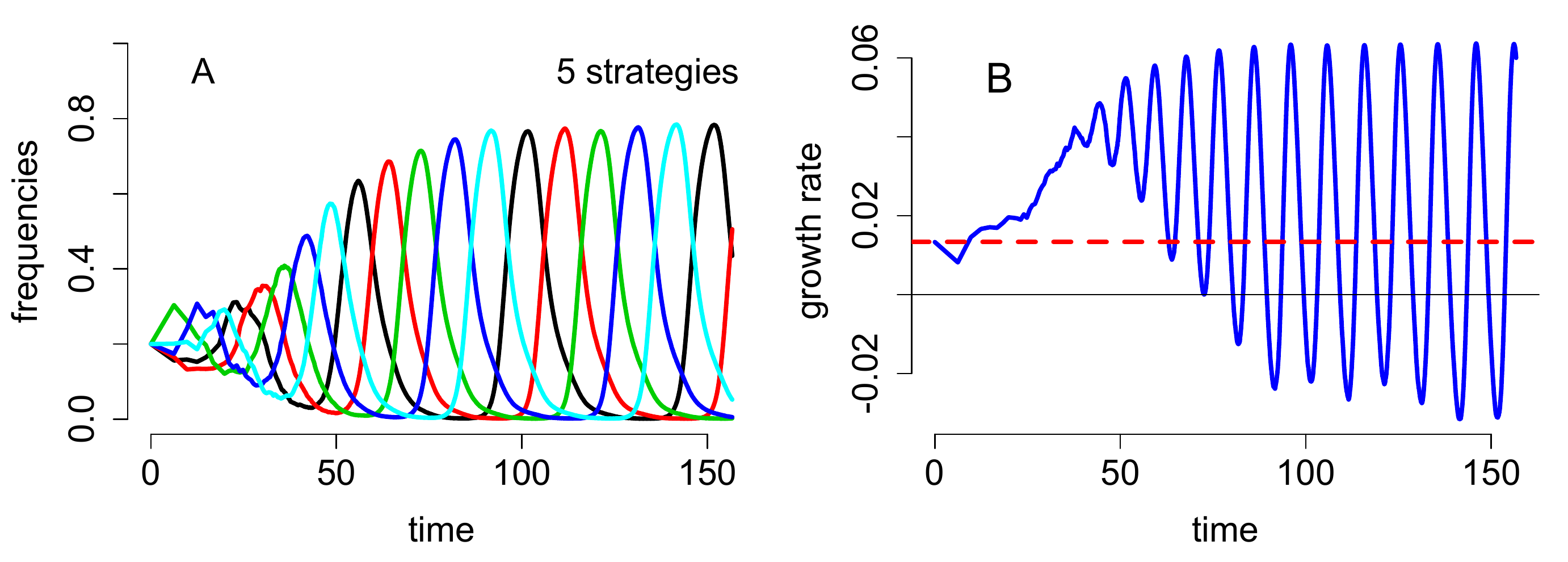}
\caption{Sample trajectories of the replicator process $x(n)$ with $b=1$, $d=2.5$, $\nu=4$, $n=5$, $B$ and $D$ given by equations (\ref{eq:example}). Since $(b-d+2\nu/5)/(b+d+\nu)=1/75>0$, there is growth and convergence with positive probability to an equilibrium attractor a non-equilibrium attractor for $n=5$ (A). In (B),  oscillations of the growth function along the stochastic trajectory are shown. The dashed line corresponds to the equilibrium value $1/75$ of the growth function. }\label{fig:replicator}
\end{figure}

To illustrate the use of this result, we consider the hypercycle replicator equations featured in \citet[Chapter 12]{hofbauer-sigmund-98}. In this game of $k$ strategies, interactions between  an individual playing strategy $i$ with individuals playing strategy $i-1$ ($n$ in the case $i=1$) catalyze births. These interactions occur at rate $\nu$, while births and deaths independent of interactions occur at rates $b$ and $d$, respectively. The birth and death matrices associated with interactions are given by
\begin{equation}\label{eq:example}
B=\begin{pmatrix}
0&0&0& . & . & . & 1\\
1&0&0& . & . & . & 0\\
0&1&0& . & . & . & 0\\
.&.&.& . & . & . & .\\
0&0&0& . & . & 1 & 0\\
	\end{pmatrix} \mbox{ and } D=0.
\end{equation}
For the replicator equation \eqref{eq:replicator} with $A=2\nu\gamma(B-D)$, \citet[Theorem 12.1.2]{hofbauer-sigmund-98} prove that there is a globally stable equilibrium $\hat x = (1/k,1/k,\dots,1/k)$ in $S_k \setminus\partial S_k$ for $k\le 4$. In contrast, for $k\ge 5$, this equilibrium is unstable (see, e.g., \citep[Section 12.1]{hofbauer-sigmund-98} and there is  is a global, non-equilibrium attractor $\A\subset S_k\setminus \partial S_k$ which attracts all initial conditions in $S_k \setminus \partial S_k$ except those on the stable manifold of $\hat x$  \citep[Theorem 12.3.1]{hofbauer-sigmund-98}.

For our stochastic analog of this replicator dynamic, we get growth and convergence with positive probability to one of the vertices whenever $b>d$. When 
$b-d+2\nu/n>0$, we get growth and convergence with positive probability to $\hat x$ whenever $n\le 4$ and growth and convergence with positive probability to the non-equilibrium attractor $\mathcal A$ whenever $n\ge 5$. Moreover, if $n\ge 5$, $b-d<0$ and $b-d+2\nu/n>0$ is sufficiently close to zero, then we conjecture that there are points in $\A$ where the growth function $f(x)=\frac{b-d}{b+d+\nu} +  \frac{2\nu}{b+d+\nu}\left(x_1x_n+x_2x_1+\dots+ x_{n-1}x_{n-2}+x_nx_{n-1}\right)$ is negative. In this case, our Theorem~\ref{thm:attractor} implies convergence with positive probability while Theorem~\ref{th:BST} from the earlier work of \citet{BenSchTar04} would not.  Figure~\ref{fig:replicator} illustrates converge to the non-equilibrium attractor for $n=5$ and the oscillations exhibited in the growth function $f(x)$.

\subsection{Selection-mutation processes} Two key evolutionary processes are natural selection in which there is differential survival or reproduction amongst genetically distinct individuals and mutation in which parents produce offspring with novel genotypes from their own. A simple, continuous time deterministic model of these processes acting simultaneously on populations tracks the gametes of the populations and assumes that selection acting on diploid individuals and mutation acting on gametes occur independently of one another. While this simplification is somewhat unrealistic, it turns out to be sufficiently realistic to provide useful insights and, for our purposes, illustrates how urn models can be used in the context of population genetics. 

We consider a population of haploid individuals (gametes) with a single locus with $k$ possible alleles, $A_1,\dots, A_k$. All gametes die at a constant rate $d>0$ and fuse with another gamete at a constant rate $\nu>0$. When two gametes, say of types $A_i$ and $A_j$, fuse to produce an individual with genotype $A_i A_j$ they produce, on average, $f_{ij}/2$ gametes of type $A_i$ and $f_{ij}/2$ gametes of  type $A_j$. For each pair $i,j$, let $B_{ij}(1),B_{ij}(2),\dots$ be a sequence of independent, identically distributed random variables taking values in $\{0,2,\dots,2m\}$ and satisfying $\mathbb{E}[B_{ij}(n)]=f_{ij}/2$. At the $n$-th fusion of a gamete of type $i$ and $j$,  $B_{ij}(n)$ gametes of type $i$ and $j$ are added to the population. Gametes of type $i$ mutate to type $j$ at rate $\mu_{ij}>0$. Namely, on this event, a gamete of type $i$ is replaced with a gamete of type $j$. Let $\mu = \sum_{i \neq j} \mu_{ij}$ be the total mutation rate.

As in the case of the replicator processes, we are interested in the discrete-time embedded stochastic process where $z(n)$ is the state of the population immediately following the $n$-th demographic event. This process has three types of demographic events: a gamete is removed due to death ($-e_i$), a gamete changes types due to mutation ($e_i-e_j$ for $i\neq j$), or new gametes are added to the population due to births following the fusing of two gametes ($l(e_i+e_j)$ for some $l \leq m$). Setting $\gamma=\frac{1}{d+\mu+\nu}$, the limiting transition functions $p_w$, as $|z|\to\infty$, corresponding to these transitions are defined by  \[
\begin{aligned}
p_{e_j-e_i}(x)&=\gamma \mu_ {ij} x_i\mbox{ for }i\neq j  \\
p_{-e_i}(x)&=\gamma d x_i\\
 p_{\ell (e_i+e_j)}(x) &= \gamma \nu x_i x_j \mathbb{P}[B_{ij}(n)=2\ell] \mbox{ for }\ell \in \{0,1,\dots,m\} 
 \end{aligned}
\]
The transition functions, $\Pi(z,w)$, for the actual process $z(n)$ are given by $p_w(x=\frac{z}{|z|},w)$ with the exception that the terms $x_i^2=\frac{z_i^2}{|z_i|^2}$ are replaced by $x_i (x_i|z|-1)/|z|$ to ensure that interactions between different individuals playing the same strategy $i$ also take place. 

Let $F=(f_{ij})_{i,j}$ and $M=(\mu_{ij})_{i,j}$. Then the mean limit ODE for $z(n)$ is given by the mutation-selection equation~\citep[Section 20.1] {hofbauer-sigmund-98}
\begin{equation}\label{eq:sm}
\frac{dx}{dt} = \gamma \nu x\circ (F x -x^T Fx) + \gamma \mu (M^Tx-x)
\end{equation}
\cite{hofbauer-85}  proved that the selection-mutation equation for \eqref{eq:sm} can exhibit gradient-like dynamics or non-equlilibrium dynamics depending on the mutation rates. For example, the following result shows that if rate of mutating to gamete type $j$ is the same for all gamete types, then the dynamics are gradient-like. 

\begin{theorem}[Hofbauer 1985] If $\mu_{ij}=\mu_j$ for all $i\neq j$, then all solutions $x(t)$ of \eqref{eq:sm} converge to the set of equilibria of \eqref{eq:sm}.\end{theorem}

Using the earlier work, Theorems~\ref{th:BST} and \ref{th:BST2} of \citet{BenSchTar04}, we get the following corollary:

\begin{corollary} Assume that $\mu_{ij}=\mu_j$ for all $i\neq j$ and $F, M$ are such that  \eqref{eq:sm} has a finite number of stable, hyperbolic equilibria $\hat x_1, \dots, \hat x_m$ of which $\hat x_1,\dots, \hat x_s$ with $s\le m$ are linearly stable. The fertility selection process $x(n)$ grows and converges to $\hat x_i$ for some $1\le i\le s$ with positive probability if 
\[
\nu\hat x_i^T F \hat x_i >d.
\]
Furthermore, on the event of linear growth, $x(n)$ converges to $\hat x_i$ for some $1\le i \le s$. 
\end{corollary}

\begin{figure}
\begin{center}\includegraphics[width=0.75\textwidth]{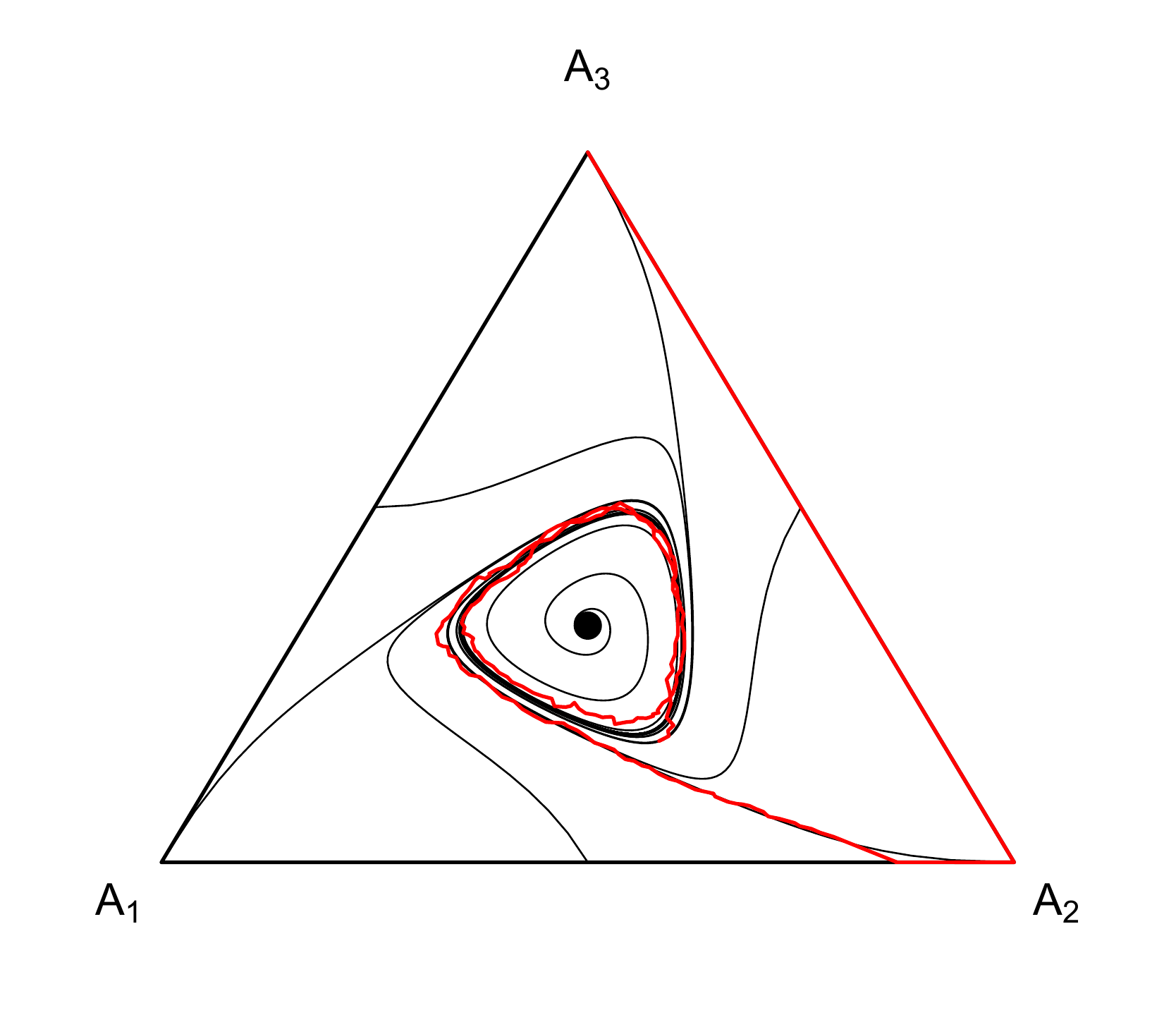}\end{center}
\caption{The dynamics of the selection-mutation process resulting in a limit cycle as described in the text. Trajectories of the mean limit ODE shown in black. A stochastic trajectory is shown in red.}
\end{figure}

Without the strong assumption on mutation rates, the selection-mutation dynamics can give rise to non-equilibrium dynamics. When this occurs, the following result is useful.

\begin{corollary}\label{cor:sm} Assume \eqref{eq:sm} has a stable periodic attractor $\A \subset S_k \setminus \partial S_k$ of period $T$. Let $x.t$ denote a solution of \eqref{eq:sm} with $x\in \A$. If
\begin{equation}\label{eq:sm2}
\frac{\nu}{T} \int_0^T (x.t)^T F (x.t)  >d 
\end{equation}
and $x \in S_k\setminus \partial S_k$, then 
\[\P \left[\sum_{n} \frac{1}{|z(n)|^{1+ \delta}} < + \infty \, \,  \forall \delta>0, \; \mbox{ and } \,  L(x(n)) \subset \mathcal{A} \right]>0.\]
Alternatively, if inequality \eqref{eq:sm2} is  reversed, then 
\[\mathbb{P} \left[ |z(n)| \geq M \; \mbox{ for n large enough} \, \mbox{ and } \;  L(x(n)) \subset \A  \right]=0.\]
\end{corollary} 

\citet{hofbauer-85} illustrated how selection-mutation equations for $k=3$ alleles can lead to oscillatory dynamics. Specifically, assume all heterozygous $A_iA_j$ with $i\neq j$ have the same fitness, i.e. $f_{ij}=f$ for all $i\neq j$,  all homozygotes $A_iA_i$ have the same fitness, i.e. $f_{ii}=f+s$ for all $i$, and the mutation rates are cyclic symmetric, i.e. $\mu_{ij}=\mu_{i-j}$ where $\sum_{i=0}^2 \mu_i =1$ and 
 the index $i$ is considered as a residue modulo $3$.  If $\mu_1\neq \mu_2$ and  $s$ is slightly larger than $\frac{9}{2}(\mu_1+\mu_2)$, then there is a stable periodic orbit (of  say period $T$). 
Assume $x$ lies on this periodic orbit. If   
\[
\frac{\nu}{T}\int_0^T s \|x.t\|^2 dt + \nu f >d
\]
then Corollary~\ref{cor:sm} implies that there is a positive probability the population grows and its distribution converges to this periodic orbit. Conversely, if 
\[
\frac{\nu}{T}\int_0^T s \|x.t\|^2 dt + \nu f <d
\]
then the population distribution can not converge to this periodic orbit.

\section{Appendix} 
In this Appendix, we prove some of the key technical lemmas and estimates used in the proofs of our main results.

\subsection{Proof of Lemma \ref{lm:diff}}  On the event $\left\{|z(r_k)| \geq \frac{4BCL}{\delta} \right\}$, we have $\Gamma_2(k,T_0) \leq \frac{2}{B^{-1}|z(r_k)|} \leq \delta/2C$. By (\ref{eq:erg}), the left expression in (\ref{eq:upbound}) is smaller than
\[\mathbb{P} \left[\Gamma_1(k,T_0) \geq \frac{\delta}{2C} \Big|  \mathcal{F}_{r_k} \right].\]
On the other hand, we have
\begin{eqnarray*}
\mathbb{E} \left[\Gamma_1(k,T_0)^2 \Big| \, \mathcal{F}_{r_k} \right] &\leq&  \mathbb{E} \left[\sum_{r=r_k}^{r_{k+1}-1} \frac{\|U\|^2}{|z(r)|^2} \right]\\
&\leq&  \frac{\|U\|^2}{B^{-1}|z(r_k)|} \mathbb{E} \left[ \sum_{r=r_k}^{r_{k+1}-1} \frac{1}{|z(r)|}\right]\\
&\leq& \frac{BT_0 \|U\|^2}{|z(r_k)|},
\end{eqnarray*}
and (\ref{eq:upbound}) holds. \qed

\subsection{Proof of Proposition \ref{pr:tech}}  By the Markov property, it suffices to prove the estimate for $k=0$ with $x(0)=x \in U$. Define $M_1=\max\{ 3B(L+1/T_0)/\delta, 12 L^2BC/\delta \}$. We have
\begin{eqnarray*}
&& \frac{1}{T_0}\left| \int_0^{T_0} f(x.s) ds  - \sum_{r=0}^{m(T_0)-1} \frac{1}{|z(r)|} f(x.\tau(r)) \right|\\
&\leq& \frac{1}{T_0} \left| \sum_{r=0}^{\text{\tiny{$\tau(r+1)$}}} \left(\int_{\tau(r)}^{\text{\tiny{$\tau(r+1)$}}} f(x.s) ds  -  \frac{1}{|z(r)|} f(x.\tau(r))\right) \right|+  \frac{1}{T_0}\left| \int_{T_0}^{\text{\tiny{$\tau(m( T_0))$}} } \text{\small{$f(x.s)ds$}} \right|\\
&\leq& \frac{1}{T_0} \left| \sum_{r=0}^{m(T_0)-1} \left(\int_{\tau(r) }^{\tau(r+1) } f(x.s) ds  -  \frac{1}{|z(r)|} f(x.\tau(r))\right) \right|+ \frac{1}{T_0} (\tau(m( T_0)) - T_0)\\
&\leq& \frac{1}{T_0}  \sum_{r=0}^{m(T_0)-1} \int_{\tau(r)}^{\tau(r+1) } \left\|f(x.s) ds  -   f(x.\tau(r)) \right\|ds+\frac{1}{T_0} \frac{1}{|z(m(T_0))|}\\
\end{eqnarray*}
\begin{eqnarray*}
&\leq& \frac{1}{T_0}  \sum_{r=0}^{m(T_0)-1} \int_{\tau(r)}^{\tau(r+1)} L|\tau(r+1) - \tau(r)| ds+ \frac{1}{T_0}\frac{1}{|z(m(T_0))|}\\
&\leq&  \frac{1}{T_0}  \sum_{r=0}^{m(T_0)-1}  \frac{L}{|z(r+1)|^2} + \frac{1}{T_0} \frac{1}{|z(m(T_0))|}\\
&\leq&   \frac{1}{|z(0)|} \left(B(L  + 1/T_0) \right)\\
& \leq & \frac{\delta}{3}
\end{eqnarray*}
on the event $V_0$.

Since
\[
\frac{1}{T_0} \left| \sum_{r=0}^{m(T_0)-1} \frac{1}{|z(r)|}\left( f(x.\tau(r))-f(x(r))\right)\right| \le L \sup_{0\le r \le m(T_0)-1} \| x.\tau(r)-x(r)\|,
\]
Lemma 2 implies 
\[
\mathbb{P}\left[\frac{1}{T_0} \left| \sum_{r=0}^{m(T_0)-1} \frac{1}{|z(r)|}\left( f(x.\tau(r))-f(x(r))\right)\right| \ge \delta/3\right]
\le \frac{9 L^2C_0(T_0)}{|z(0)| \delta^2} 
\]
on the event $V_0$. 

Finally, 
\[
\frac{1}{T_0}\left| \int_0^{T_0} f(x.s)-f(y.s) ds \right|  \le L \sup_{0\le s \le T_0} \|x.s-y.s\| \le \frac{\delta}{3}
\]
for $x\in U$ by choosing $U$ to be a sufficiently small neighborhood of $y$. 

These three estimates plus the triangle inequality complete the proof of the proposition, with $C_1(T_0) = 9L^2C_0(T_0)$. \qed

\subsection{Proof of Proposition \ref{pr:prop2}}
By the Markov property, it suffices to prove the estimate for $k=0$ with $x(0)=x \in U$. Define  $N(i)=|z(i+1)|-|z(i)|$, $D(i)=N(i)-E[N(i) | z(i)]$ and
$$G= \frac{1}{T_0} \sum_{i=0}^{m(T_0)-1} \frac{D(i)}{|z(i)|}.$$
Observe that $|D(i)|\leq 2m$. Therefore
\begin{eqnarray*}
E\left[G^2 | z(0) \right] &&\le \frac{1}{T_0^2}
E \left[ \sum_{i=0}^{m(T_0)-1} \frac{D(i)^2}{|z(i)|^2} \mid z(0)\right] \\
&& \le\frac{4m^2B^2}{T_0^2|z(0)|^2} m(T_0) \leq \frac{4m^2B^3}{T_0|z(0)|}
\end{eqnarray*}
where we have used the fact that $T_0 B^{-1} |z(0)|\leq
m(T_0)\leq T_0 B |z(0)|$ (see point $(ii)$ before Proposition \ref{pr:prop2}). Chebyshev's inequality implies 
\begin{equation}\label{eq:cheby}
P \left[|G|\geq \frac{\delta}{3} \mid z(0) \right] \leq \frac{36m^2B^3}{T_0|z(0)|\delta^2}.
\end{equation}

Let $C_1(T_0)$ and $M_1$ be as defined by Proposition~\ref{pr:tech}. Notice that, by definition of $f$ and assumption (A2), we have
\[\left|\mathbb{E}(N(i) \mid z(i)) - f(x(i)) \right| \leq \frac{a}{|z(i)|}.\]
Hence, 
\[\frac{1}{T_0}\left|  \sum_{i=0}^{m(T_0)-1} \frac{\mathbb{E}(N(i) \mid z(i)) - f(x(i))}{|z(i)|} \right| \leq \frac{aBT_0}{|z(0)|} \leq \frac{\delta}{3}\]
provided $M_2$ is large enough.

Consequently, Proposition~\ref{pr:tech} with a $\delta$ value of $\delta/3$  and inequality \eqref{eq:cheby} imply 
\[
\begin{aligned}
& \mathbb{P}\left[ \frac{1}{T_0}\left|  \sum_{i=0}^{m(T_0)-1} \frac{N(i)}{|z(i)|}-\int_0^{T_0} f(x(0).s)ds\right| \ge \delta \mid z(0)  \right]\\
&\le  \mathbb{P}\left[ \frac{1}{T_0}\left|  \sum_{i=0}^{m(T_0)-1} \frac{D(i)}{|z(i)|}\right| \ge \frac{\delta}{3} \mid z(0)\right]\\
&+ \mathbb{P}\left[\frac{1}{T_0}\left|  \sum_{i=0}^{m(T_0)-1} \frac{f(x(i))}{|z(i)|}-\int_0^{T_0} f(x(0).s)ds\right| \ge \frac{\delta}{3} \mid z(0)\right]\\
&\le \frac{36m^2B^3}{T_0|z(0)|\delta^2}+\frac{9C_1(T_0)}{|z(0)|\delta^2} \le \frac{C_2(T_0)}{|z(0)|\delta^2}
\end{aligned}
\]
on the event $W_0$ with $M_2\geq M_1$ large enough  and $C_2(T_0)=9C_1(T_0)+ \frac{36m^2B^3}{T_0}$. \qed

\subsection{Proof of claim (\ref{claim})} Let us prove (\ref{claim}): 
\begin{equation}  \label{claim2} \sum_{i=r_k+1}^{r_{k+1}-1} \frac{N(i)}{|z(i)|} \leq  \frac{n_k}{B^{-1}|z(r_k)|} +  \frac{m}{2}(r_{k+1} - r_k) \left(\frac{1}{ B^{-1}|z(r_k)|} - \frac{1}{B^{-1}|z(r_k)|+m}\right)
\end{equation}
where $n_k = |z(r_{k+1})| - |z(r_k)|$. Assume without loss of generality that $n_k >0$.  For every ball which is removed between times $\tau(r_k)$ and $\tau(r_{k+1})$ when state is $z(r)$, a ball is added when the state is at least $|z(r)|-m$.  The quantity
\[\frac{1}{|z(r)|-m} - \frac{1}{|z(r)|}\]
is maximal when $|z(r)|-m$ is minimal, and $z(r)-m$ can never be smaller than $B^{-1}|z(r_k)|$ for $r \in [r_k,r_{k+1}-1]$.  Moreover there is at most $m(r_{k+1}-r_k)/2$ balls removed in the process). Finally the remaining $n_k$ balls were added when the state was at least $B^{-1}|z(r_k)|$. A very rough upper bound is given by (\ref{claim2}). $\; \; \square$

\paragraph{Acknowledgements} This work was supported in part by U.S. National Science Foundation Grants EF-0928987 and DMS-1022639 to Sebastian Schreiber. The authors would also like to thank the Aix-Marseille School of Economics for providing travel funds to Mathieu Faure to work on this project in S. Schreiber's lab.

\section*{References}

\bibliographystyle{elsarticle-harv} 

\bibliography{qsd_biblio}

\end{document}